\documentclass{article}

\usepackage{amsmath,amssymb,indentfirst,alltt,theorem}
\newcommand{\R}{{\mathbb{R}}}
\newcommand{\N}{{\mathbb{N}}}
\def\bgneqy{\begin{eqnarray}}
\def\endeqy{\end{eqnarray}}
\def\bgneqy*{\begin{eqnarray*}}
\def\endeqy*{\end{eqnarray*}}
\newcommand{\C}{{\mathbb{C}}}

\newtheorem{thm}{Theorem}
\newtheorem{Th}[thm]{\normalsize{\bf{Theorem}}}

\newtheorem{Cor}[thm]{\normalsize{\bf{Corollary}}}
\newtheorem{Lem}[thm]{\normalsize{\bf{Lemma}}}
\newtheorem{Prop}[thm]{\normalsize{\bf{Proposition}}}
\newtheorem{Def}[thm]{\normalsize{\bf{ Definition}}}

\newtheorem{Ex}[thm]{Example}

\newcommand\dint{\displaystyle\int}
\newcommand\dlim{\displaystyle\lim}
\newcommand\dsup{\displaystyle\sup}

\RequirePackage[colorlinks,hyperindex]{hyperref}
\begin{document}
%\articletype{}
\title{\textbf{Sobolev type spaces associated with the Poly-axially operator}}
\author{{Belgacem Selmi \textsuperscript{a} \thanks{ Email:
Belgacem.Selmi@fsb.rnu.tn},\quad Chahiba Khelifi \textsuperscript{b}
\thanks{ Email: khelifi\_chahiba@yahoo.com}}
\date{ \small{
\textsuperscript{a} Faculty of Sciences of Bizerte, Carthage University, Tunisia;\\
\textsuperscript{b} Faculty of Sciences of Bizerte, Carthage
University, Tunisia.}}}
%\date{}

\maketitle

\begin{abstract}
In this paper n-dimensional Sobolev type spaces $
E_{\alpha}^{s,p}(\R^n_+)$ $(\alpha\in \R^n,\;\;\alpha_1>
-\frac{1}{2},...,\alpha_n>-\frac{1}{2}, s\in \R, p\in [1,+\infty])$
are defined on $\R^n_+$ by using Fourier-Bessel transform. Some
properties including completeness and embedding results for these
spaces are obtained, Poincar\'{e}'s inequality and Reillich theorem
are proved.
\end{abstract}

\section{Introduction}

The important role of Sobolev type spaces in analysis and its
application is well known. In the classical case of $\R$, R. Adams
\cite{ad} introduced the Sobolev type spaces $E_*^{s,p}(\R)$ given
by
$$E_*^{s,p}(\R)=\left\{T\in {\cal S}'_*(\R)| (1+\xi^2)^{s}{\cal F}_{0}(T)\in L^p(\R) \right\},\quad s\in\R,\;\;1\leq p\leq+\infty,$$
where ${\cal F}_0$ is the even Fourier transform given by ${\cal
F}_0(f)(\lambda)=\dint_0^{+\infty}f(x)\cos(\lambda x)dx$  and ${\cal
S}'_*(\R)$ is the space of tempered and even distributions.\\
This theory of Fourier transform on the space of tempered
distributions on $\R^n$ introduced  by L. Schwartz in \cite{sch} has
been exploited by many authors in the study of Sobolev spaces
\cite{ad,hb,ls}. The same theory has been generalized for other
spaces and it is applied to
investigate Sobolev type spaces \cite{Assal,B.Sal,pp.s,jel}.\\
In this paper, we consider the n-orders differential Bessel operator
defined for $\alpha=(\alpha_1,...,\alpha_n)\in \R^n$ such that
$\alpha_i>-\frac{1}{2}$ for $i=1,...,n$, by
$$\Delta_{\alpha}= \sum_{i=1}^{n}\frac{\partial^2}{\partial x_i^2}+ \frac{2\alpha_i+1}{x_i}\frac{\partial}{\partial
x_i}.$$
 which is called the Poly-axially operator (see \cite{H}) and for the two-dimensional case one can see \cite{ab.k,zf}.
The main object of this work is to generalize the subject of the
Sobolev type spaces $E_{\alpha}^{s,p}(\R^n_+)$ to the case of
harmonic analysis associated with the Poly-axially operator
$\Delta_{\alpha}$. Then, we define the Sobolev type space
$E_{\alpha}^{s,p}(\R^n_+)$, $(s\in \R,\; p\geq 1)$, by using the
Fourier-Bessel transform ${\cal{F}}_{\alpha}$ \cite{kt} given by
$${\cal{F}}_{\alpha}f(\lambda)= \dint_{\R^n_+} f(x_1,...,x_n)
j_{\alpha}(\lambda_1 x_1,...,\lambda_n
x_n)d\mu_{\alpha}(x_1,...,x_n),$$ where
 $$d\mu_{\alpha}(x)=\prod_{i=1}^{n}x_i^{2\alpha_i+1}dx_i\;\; and\;\;\;j_{\alpha}(\lambda_1 x_1,...,\lambda_n
x_n)=\prod_{i=1}^{n}j_{\alpha_i}(\lambda_ix_i),\;\;for
\;\;x,\lambda\in \R^n_+,$$ \Big( $j_{\gamma},\;\gamma>-\frac{1}{2}$,
is the normalized Bessel function of first kind and order $\gamma$
\cite{wa}\Big). Moreover, we investigate its properties and some
results of the Sobolev type spaces. Next, in the case $p=2$, we
denote by $H_{\alpha}^{s}(\R^n_+)$ the Sobolev type spaces and we
study its properties. Finally, we prove Reillich theorem and
Poincar\'{e} inequality for these spaces.

The paper is organized as follows: in the first section we recall
some results in harmonic analysis which we need in the sequel. The
second section contains the definition of Sobolev type spaces
$E_{\alpha}^{s,p}(\R^n_+)$ associated with the Poly-axially operator
$\Delta_{\alpha}$ and some of their properties. In the third section
we study in particular the space $E_{\alpha}^{s,2}(\R^n_+)$ denoted
by $H_{\alpha}^s(\R^n_+)$ then Reillich theorem and Poincar\'{e}'s
inequality and we give some applications.

%%%%%%%%%%%%%%%%%%%%%%%%%%%%%%%%%%%%%%%%%%%%%%%%%%%%%%%%%%%%%%%%%%%%%%%%%%%%%%%%%%%%%%%%%%%%%%%%%%%%%%%%%%%%%%%%%%%%%%%%%%%%%%%%%%%%%%%
%%%%%%%%%%%%%%%%%%%%%%%%%%%%%%%%%%%%%%%%%%%%%%%%%%%%%%%%%%%%%%%%%%%%%%%%%%%%%%%%%%%%%%%%%%%%%%%%%%%%%%%%%%%%%%%%%%%%%%%%%%%%%%%%
\section{Preliminaries}
\label{sec:2}
%%%%%%%%%%%%%%%%%%%%%%%%%%%%%%%%%%%%%%%%%%%%%%%%%%%%%%%%%%%%%%%%%%%%%%%%%%%%%%%%%%%%%%%%%%%%%%%%%%%%%%%%%%%%%%%%%%%%%%%%%%%%%%%%%

{\bf{Notations:}} We first begin by some notations of functional
spaces which we need in this work. A functional $f$ on $\R^n$ is
said to be even if it is even in each variable \cite{zf}.\\
$\bullet\;\;  \R^n_+=\left\{ x=(x_1,...,x_n)\in \R^n,\;\;
x_1>0,x_2>0,...,x_n>0\right\}$.\\
$\bullet\;\; {\cal C}_e(\R^n)$ consists of all
even continuous functions on ${\R}^n$.\\
$\bullet\;\; {\cal C}^m_e(\R^n)$ the space of all even ${\cal
C}^m$-functions on $\R^n$ and ${\cal C}_{e,0}(\R^n)$ the space of
even continuous functions defined on $\R^n$ satisfying:
$$\dlim_{\|x\|\longrightarrow +\infty}f(x)=0\;\;
and\;\;\|f\|_{{\cal C}_{e,0}}=\dsup_{x\in\R^n_+}|f(x)|<+\infty.$$
$\bullet\;\;  {\cal S}_e(\R^n)$: the Schwartz space consists of all
even ${\cal C}^\infty$-functions on $\R^n$
 which are rapidly decreasing as their derivatives, provided with the topology defined by for all $\gamma=(\gamma_1,...,\gamma_n)\in
\N^n$\\
$\rho_{k}(f)=\dsup_{x\in\R^n_+, |\gamma|\leq k}
(1+\|x\|^2)^k\Big|\partial^{\gamma}f(x)\Big|<\infty$, for all
$k \in \N$,\\
with
$\partial^{\gamma}f(x)=\dfrac{\partial^{|\gamma|}}{\partial^{\gamma_1}_{x_1}...\partial^{\gamma_n}_{x_n}}f(x)$,
$|\gamma|=\gamma_1+...+\gamma_n$ and we denote by
$\|x\|^2=x_1^2+...+x_n^2$.\\
$\bullet\;\;  {\cal S}'_e(\R^n)$ the space of even tempered
distributions on $\R^n$.\\
$\bullet\;\;  L_{ \alpha}^p(\R^n_+)$, for
$\alpha=(\alpha_1,...,\alpha_n)\in \R^n$ such that
$\alpha_1\geq-\frac{1}{2},...,\alpha_n\geq-\frac{1}{2}$
 and $1\leq p<\infty$ (resp. $p=\infty$), the space of measurable functions
 $f$ on $\R^n_+$ such that,
$$\|f\|^p_{L_{\alpha }^p}
=\dint_{\R^n_+}|f(x)|^pd\mu_{\alpha}(x)<+\infty,\;\;\Big(resp.\;\;
\|f\|_{\infty}=\|f\|_{L_{ \alpha
}^\infty}=ess\dsup_{x\in\R^n_+}|f(x)|<+\infty \Big),$$ where
$$d\mu_{\alpha}(x)=x_1^{2\alpha_1+1}...x_n^{2\alpha_n+1}dx_1...dx_n.$$
$\bullet\;\;  {\cal D}_e(\R^n)$ the space of all even ${\cal
C}^\infty$-functions with compact support.\\

For $x=(x_1,...,x_n)\in \R^n_+$, we have
\begin{equation}
\Delta_{\alpha}[\prod_{i=1}^{n}j_{\alpha_i}(\xi_i
x_i)]=-\|\xi\|^2\prod_{i=1}^{n}j_{\alpha_i}(x_i)
\end{equation}

One can remark that\\
$\bullet\;\;$ If $\alpha_i=-\frac{1}{2}$ for $i=1,...,n$ then
$$\Delta_{\alpha}= \sum_{i=1}^{n}\frac{\partial^2}{\partial
x_i^2}\;\;\;is\; the\; Laplace\; operator\; on\; \R^n.$$
 and for all $x\in \R^n_+$ we have\\
$$j_{-\frac{1}{2},...,-\frac{1}{2}}(x_1,...,x_n)=\prod_{i=1}^{n}j_{-\frac{1}{2}}(x_i)=\prod_{i=1}^{n}\cos(x_i).$$
$\bullet\;\;$ If $\alpha_i=-\frac{1}{2}$ for $i=1,...,n-1$ and
$\alpha_n>-\frac{1}{2}$, then
$$\Delta_{\alpha}= \sum_{i=1}^{n}\frac{\partial^2}{\partial x_i^2}+\frac{2\alpha_n+1}{x_n}\frac{\partial}{\partial x_n},$$
 $\Delta_{\alpha}$ coincides with the Weinstein operator defined on $\R^n_+$ \cite{Najib,Garna}.

\subsection{Generalized Bessel translation operator}
The generalized translation operator associated with the
Poly-axially operator $\Delta_{\alpha}$ is defined \cite{H}: for all
$x=(x_1,...,x_n),y=(y_1,...,y_n)\in\R^n_+$ by

$${\cal T}^{\alpha}_yf(x)=c'_{ \alpha}\dint_{[0,\pi]^n}
f(X_1,...,X_n)\prod_{i=1}^{n}(\sin\theta_i)^{2\alpha_i}d\theta_1...d\theta_n,$$
with $c'_{
\alpha}=\prod_{i=1}^{n}\frac{\Gamma(\alpha_i+1)}{\pi^{n/2}\Gamma(\alpha_i+1/2)}$
and $X_i=\sqrt{x_i^2+y_i^2-2y_ix_i\cos\theta_i},$ for $i=1,...,n$.
\\
Using the following change of variables: $ z_1=X_1,...,z_n=X_n$, the
generalized translation operator can be rewritten as: for $x,y\in
\R^n_+$
\begin{equation}
{\cal T}^{\alpha}_yf(x)=\dint_{\R^n_+} w_{ \alpha}(x,y,z)f(z)d\mu_{
\alpha}(z),
\end{equation}
where  the kernel $w_{\alpha}$ is given by\\
$$w_{\alpha}(x,y,z)=
 \left\{\begin{array}{c}
\frac{c'_{\alpha}}{2^{2|\alpha|-n}}\frac{\prod_{i=1}^{n}
\Big(\Big[z_i^2-(x_i-y_i)^2\Big]\Big[(x_i+y_i)^2-z_i^2\Big] \Big)
  ^{\alpha_i-1/2}}{\prod_{i=1}^{n}(z_iy_ix_i)^{2\alpha_i}},\;\;
                       if\;(z_1,...,z_n)\in \prod_{i=1}^{n}[|x_i-y_i|,x_i+y_i],\\
                       \\
                        0,\;\;\;\;\;\quad\quad otherwise.
\end{array}\right.$$\\

For $x,y,z \in \R^n_+$, the kernel $w_{\alpha}$ have the following
properties (for the univariate case we can see \cite{kt}).
\begin{Prop} For $\alpha_i>-\frac{1}{2}$ for $i=1,...,n$ , we have:
\begin{enumerate}
\item $\forall x,z\in \R^n_+$;  $w_{\alpha}(x,z,u)\geq 0,\;\;\forall u\in \R^n_+.$
\item $\forall x,z\in \R^n_+$;
$supp(w_{\alpha}(x,z,\cdot))\subset
\prod_{i=1}^{n}[|x_i-z_i|,x_i+z_i]$.
\item $\forall x,z,u\in \R^n_+$;  $w_{\alpha}(x,z,u)=w_{\alpha}(x,u,z)=w_{\alpha}(z,x,u)$.
\item $\forall x,u\in \R^n_+$;
 $\dint_{\R^n_+}w_{\alpha}(x,u,z)d\mu_{\alpha}(z)=1.$
 \end{enumerate}
\end{Prop}
In particular the following product formula holds
$$j_{\alpha}(t_1 x_1,...,t_n x_n)j_{\alpha}(t_1 z_1,...,t_n z_n)= {\cal T}_z(j_{\alpha}(t_1 \cdot,...,t_n \cdot))(x),
\quad for\;all\; x,z\in\R^n_+.$$ If $f\in L_{\alpha}^p(\R^n_+)$ then
for all $x\in \R^n_+$, ${\cal T}_xf\in L_{\alpha}^p(\R^n_+)$ and we
have
\begin{equation}\label{e:p.1}
\|{\cal T}_xf\|_{L_{\alpha}^p}\leq \|f\|_{L_{\alpha}^p}.
\end{equation}

\subsection{ n-Dimensional Fourier-Bessel transform}
The Fourier-Bessel transform ${\cal F}_{\alpha}$ is defined for
suitable function $f$ as follows:
$${\cal{F}}_{\alpha}f(\lambda)= \dint_{\R^n_+} f(x_1,...,x_n)
j_{\alpha}(\lambda_1 x_1,...,\lambda_n
x_n)d\mu_{\alpha}(x_1,...,x_n).$$ Moreover, we prove that the
Fourier-Bessel transform ${\cal F}_{ \alpha}$ satisfies the
following properties.
\begin{Prop}\label{p1} For all $f \in L_{\alpha}^1(\R^n_+)$, we
have
\begin{enumerate}
\item ${\cal{F}}_{\alpha}(f)\in {\cal C}_{e,0}(\R^n)$ and
$\|{\cal{F}}_{\alpha}(f)\|_{\infty}\leq \|f\|_{L_{\alpha}^1}$.
\item For $x\in \R^n_+$,
\begin{equation}\label{e:e}
{\cal{F}}_{\alpha}(\Delta_{\alpha}f)(x)=-\|x\|^2{\cal{F}}_{\alpha}f(x).
\end{equation}
\end{enumerate}
\end{Prop}
\begin{Th}(Inversion formula)\label{IF}
\begin{enumerate}
\item If $f \in  L_{\alpha}^1(\R^n_+)$ such that
${\cal{F}}_{\alpha}(f) \in  L_{\alpha}^1(\R^n_+)$,
 then for almost all $x\in \R^n_+$ we have
$$f(x)=c^2_{\alpha}{\cal{F}}_{\alpha}{\cal{F}}_{\alpha}f(x),$$
where
$c_{\alpha}=\frac{2^{-|\alpha|}}{\prod_{i=1}^{n}\Gamma(\alpha_i+1)}$.
\item ${\cal{F}}_{\alpha}$ is a topological isomorphism from
${\cal S}_e(\R^n)$ into itself and
${\cal{F}}_{\alpha}^{-1}=c^2_{\alpha}{\cal{F}}_{\alpha}$.
\end{enumerate}
\end{Th}
The Fourier-Bessel transform ${\cal F}_{\alpha}$ can be extended to
$L_{\alpha}^2(\R^n_+)$ and we have
\begin{Th}(Plancherel formula)\\
The Fourier-Bessel transform  ${\cal F}_{\alpha}$ is an isomorphism
of $L_{ \alpha,\beta }^2(\R^n_+)$ and $for \;all \;\;f\in
L_{\alpha}^2(\R^n_+)$ we have
 $$\|f\|_{L_{\alpha}^2}=c_{\alpha}\|{\cal{F}}_{\alpha}f\|_{L_{\alpha}^2}.$$
\end{Th}
Now, the following proposition holds.
\begin{Prop}\label{p2} For all $f,g\in L_{\alpha}^1(\R^n_+)$, we have
\begin{enumerate}
\item $\dint_{\R^n_+} f(\xi){\cal
F}_{\alpha}(g)(\xi)d\mu_{\alpha}(\xi)= \dint_{\R^n_+} g(\xi){\cal
F}_{\alpha}(f)(\xi)d\mu_{\alpha}(\xi)$.
\item ${\cal F}_{\alpha}({\cal T}_xf)(y)=\prod_{i=1}^{n}
j_{\alpha_i}(x_iy_i){\cal F}_{\alpha}f(y),$ for all $x,y\in\R^n_+$.
\end{enumerate}
\end{Prop}

{\bf{Remark.}} As in the one-dimensional case \cite{kt}, one can
prove that for all $p\geq 1$,  $L_{\alpha}^p(\R^n_+)$ is dense in
${\cal S}'_e(\R^n)$. Hence, for all $\psi \in L_{\alpha}^p(\R^n_+)$
and $f\in {\cal S}_e(\R^n)$, $\langle\psi,f\rangle$ means the value
of $\psi \in {\cal S}'_e(\R^n)$ on $f$ and it is given by
$$\langle\psi,f\rangle=\dint_{\R^n_+}\psi(x)f(x)d\mu_{\alpha}(x).$$

We define the n-dimensional Fourier-Bessel transform for $T\in {\cal
S}'_e(\R^n)$ by:
$$\forall \;\;\varphi\in {\cal S}_e(\R^n),\;\;\langle
{\cal{F}}_{\alpha}(T),\varphi\rangle=\langle
T,{\cal{F}}_{\alpha}(\varphi)\rangle.$$ and ${\cal{F}}_{\alpha}$ is
an isomorphism from ${\cal S}'_e(\R^n)$ into itself.

\subsection{The generalized Bessel convolution product}
The generalized Bessel convolution product is defined for suitable
functions $f$ and $g$ as follows, for all $x\in \R^n_+$:
$$f\ast_{\alpha} g(x)=\int_{\R^n_+} f(y){\cal T}_xg(y)d\mu_{\alpha}(y).$$
The generalized Bessel convolution product satisfies the following
properties
\begin{Prop}
\begin{enumerate}
\item Let $f,g \in L_{\alpha}^1(\R^n_+)$, then
$f\ast_{\alpha} g\in L_{\alpha}^1(\R^n_+)$ and we have
\begin{equation}\label{e:p.2}
{\cal F}_{\alpha}(f\ast_{\alpha} g)={\cal F}_{\alpha}(f){\cal
  F}_{\alpha}(g).
\end{equation}
\item Let $f \in L_{\alpha}^1(\R^n_+)$ such that ${\cal
F}_{\alpha}(f)\in L_{\alpha}^1(\R^n_+)$. Then for all $g\in
L_{\alpha}^1(\R^n_+)$, we have $f g\in L_{\alpha}^1(\R^n_+)$ and
\begin{equation}\label{e1}
{\cal F}_{\alpha}(fg)=c_{\alpha}^2 {\cal
F}_{\alpha}(f)\ast_{\alpha}{\cal F}_{\alpha}(g).
\end{equation}
\item For $f \in L_{\alpha}^p(\R^n_+)$, $p\in [1,+\infty)$, $g
\in L_{\alpha}^1(\R^n_+)$ we have $ f\ast_{\alpha} g\in
L_{\alpha}^p(\R^n_+)$ and
\begin{equation}\label{e.1.1}
\|f\ast_{\alpha} g\|_{L_{\alpha}^p}\leq
\|f\|_{L_{\alpha}^p}\|g\|_{L_{\alpha}^1}.
\end{equation}
\end{enumerate}
\end{Prop}

The proof of all this results follows from adaptations of those for
the associated univariate cases \cite{hb,kt,wa}.

%%%%%%%%%%%%%%%%%%%%%%%%%%%%%%%%%%%%%%%%%%%%%%%%%%%%%%%%%%%%%%%%%%%%%%%%%%%%%%%%%%%%%%%%%%%%%%%%%%%%%%%%%%%%%%%%%%%%%%%%%%%%%%%%%%%%%
%%%%%%%%%%%%%%%%%%%%%%%%%%%%%%%%%%%%%%%%%%%%%%%%%%%%%%%%%%%%%%%%%%%%%%%%%%%%%%%%%%%%%
\section{ The Sobolev type spaces $E_{\alpha}^{s,p}(\R^n_+)$}
%%%%%%%%%%%%%%%%%%%%%%%%%%%%%%%%%%%%%%%%%%%%%%%%%%%%%%%%%%%%%%%%%%%%%%%%%%%%%%%%%%%%%
%%%%%%%%%%%%%%%%%%%%%%%%%%%%%%%%%%%%%%%%%%%%%%%%%%%%%%%%%%%%%%%%%%%%%%%%%%%%%%%%%%%%%%%%%%%%%%%%%%%%%%%%%%%%%%%%%%%%%%%%%%%%%%%%%%%%%

\begin{Def}
For all $(s,p) \in \R \times [1,+\infty)$ and $\xi \in \R^n_+$, the
Sobolev type space $E_{\alpha}^{s,p}(\R^n_+)$ associated with the
operator $\Delta_{\alpha}$ is defined by:
$$E_{\alpha}^{s,p}(\R^n_+)=\left\{T\in {\cal S}'_e(\R^n)/
(1+\|\xi\|^2)^{s}{\cal F}_{\alpha}(T)\in L_{\alpha}^p(\R^n_+)
\right\}.$$
\end{Def}

\begin{Ex}
Let $(s,p) \in \R \times[1,+\infty[$ such that $2s p +
(2-p)(|\alpha| +\frac{n}{2})<-n$, then for all $x\in \R^n_+$ the
Dirac
distribution $\delta_x$ is an element of $E_{\alpha}^{s,p}(\R^n_+)$.\\

\end{Ex}
{\bf{Proof.}} Let $x\in \R^n_+$ and $\varphi\in {\cal S}_e(\R^n)$,
\begin{eqnarray*}
<{\cal F}_{\alpha}(\delta_x),\varphi>=<\delta_x,{\cal
F}_{\alpha}(\varphi)>={\cal F}_{\alpha}(\varphi)(x)
=\dint_{\R^n_+}\varphi(y)\prod_{i=1}^{n}j_{\alpha_i}(x_iy_i)d\mu_{\alpha}(y).
\end{eqnarray*}
Then, $y\longmapsto{\cal F}_{\alpha}(\delta_x)(y)$ coincides with
the bounded function $y\longmapsto
j_{\alpha}(x_1y_1,...,x_ny_n)$.\\
On the other hand, there exist a constant $C$ such that
$$|{\cal
F}_{\alpha}(\delta_x)(y)|=|j_{\alpha}(x_1y_1,...,x_ny_n)|=\prod_{i=1}^{n}|j_{\alpha}(x_iy_i)|
\leq \dfrac{C}{\prod_{i=1}^{n}(x_iy_i)^{\alpha_i+\frac{1}{2}}
},\;\;for\;\;x,y\in \R^n_+.$$

Then,
\begin{eqnarray*}
\|\delta_x\|^p_{E_{\alpha}^{s,p}}&=&\|(1+\|\xi\|^2)^s{\cal
F}_{\alpha}(\delta_x)\|^p_{L_{\alpha}^p}\\
&=&\dint_{\R^n_+}(1+\|\xi\|^2)^{sp}|{\cal
F}_{\alpha}(\delta_x)(\xi)|^pd\mu_{\alpha}(\xi)\\
&\leq&\dint_{0<\|\xi\|\leq1}(1+\|\xi\|^2)^{sp}d\mu_{\alpha}(\xi)+\dfrac{C}{\prod_{i=1}^{n}(x_iy_i)^{\alpha_i+\frac{1}{2}}
}\dint_{\|\xi\|>1}(1+\|\xi\|^2)^{sp}\Big(\prod_{i=1}^{n}\xi_i^{(2-p)(\alpha_i+\frac{1}{2})}\Big)d\xi_1...d\xi_n\\
&<& \infty.
\end{eqnarray*}

The result is proved.

\begin{Th}\label{T:1}
For any $m\in \N$ and $p\in[1,+\infty)$ we have:
$$E_{\alpha}^{m,p}(\R^n_+)=\left\{   T\in {\cal S}'_e(\R^n)\;and\;
{\cal F}_{\alpha}((-\Delta_{\alpha})^j(T))\in L_{\alpha}^p(\R^n_+)
\; for\; j\in \left\{0,...,m\right\} \right\},$$ where
$${\cal F}_{\alpha}((-\Delta_{\alpha})^j(T))(\xi)=\|\xi\|^{2j}{\cal F}_{\alpha}(T)(\xi),\;\;\xi\in \R^n_+$$
\end{Th}

{\bf{Proof.}} Let $T\in  {\cal S}'_e(\R^n)$. Then from the relation
(\ref{e:e}) we have:
$${\cal F}_{\alpha}((-\Delta_{\alpha})(T))(\xi)=
\|\xi\|^2{\cal F}_{\alpha}(T)(\xi).$$ Since
\begin{equation}\label{e:I.1}
(1+\|\xi\|^2)^m{\cal F}_{\alpha}(T)=\sum_{j=0}^{m}C_m^j{\cal
F}_{\alpha}((-\Delta_{\alpha})^j(T))
\end{equation}
and
\begin{equation}\label{e:I.2}
\forall\;j\in \left\{0,...,m\right\},\;\; |{\cal
F}_{\alpha}((-\Delta_{\alpha})^j(T))|\leq (1+\|\xi\|^2)^{m}|{\cal
F}_{\alpha}(T)|.
\end{equation}
From (\ref{e:I.1}) and (\ref{e:I.2}) we can deduce the result.

\begin{Cor}\label{c:I.1}
For any $m\in \N$, we have:
$$E_{\alpha}^{m,2}(\R^n_+)=\left\{f\in L_{\alpha}^2(\R^n_+)/ (-\Delta_{\alpha})^jf \in L_{\alpha}^2(\R^n_+),\; for\; j\in
\left\{0,...,m\right\}\right\}.$$

 In particular  $E_{\alpha}^{0,2}(\R^n_+)=L_{\alpha}^2(\R^n_+).$
\end{Cor}

{\bf{Proof.}} Theorem \ref{T:1} and the Plancherel theorem of
Fourier-Bessel transform ${\cal F}_{\alpha}$ give the result.

\begin{Prop}\label{p:I.1}
For all $(s,p)\in \R\times[1,+\infty)$. The map
\begin{eqnarray*}
\Psi:&(E_{\alpha}^{s,p},\|.\|_{E_{\alpha}^{s,p}})&\longrightarrow
 (L_{\alpha}^p(\R^n_+),\|.\|_{L_{\alpha}^p})\\
 &T&\longmapsto c_{\alpha}(1+\|\xi\|^2)^{s}{\cal
 F}_{\alpha}(T)(\xi)\\
 \end{eqnarray*}
 is an isometric isomorphism, where $\|.\|_{E_{\alpha}^{s,p}}$
 is the norm on $E_{\alpha}^{s,p}(\R^n_+)$ defined by
 $$\|T\|_{E_{\alpha}^{s,p}}=c_{\alpha}
 \|(1+\|\xi\|^2)^{s}{\cal F}_{\alpha}(T)\|_{L_{\alpha}^p}.$$
\end{Prop}

{\bf{Proof.}} We use the fact that the Fourier-Bessel transform
${\cal F}_{\alpha}$ is an isomorphism from ${\cal S}'_e(\R^n)$ into
itself.

\begin{Cor}\label{c:I.2}
For all $(s,p) \in \R\times [1,+\infty)$, the space
$E_{\alpha}^{s,p}$ endowed with the norm
 $\|.\|_{E_{\alpha}^{s,p}}$ is a Banach space.
\end{Cor}
{\bf{Remark.}} For all $(s,p) \in \R\times [1,+\infty)$, ${\cal
S}_e(\R^n)$ is dense in $E_{\alpha}^{s,p}.$

\begin{Prop}\label{p:I.2}
For all $s\in \R$, the norm $\|.\|_{E_{\alpha}^{s,2}}$ derives from
the inner product
$$(S,T)_{E_{\alpha}^{s,p}}=c_{\alpha}^2
\Big(\dint_{\R^n_+}(1+\|\xi\|^2)^{2s}{\cal
F}_{\alpha}(S)(\xi)\overline{{\cal
F}_{\alpha}(T)}(\xi)d\mu_{\alpha}(\xi)\Big).$$ In addition,
$E_{\alpha}^{s,2}$ equipped with this inner product is a Hilbert
space.
\end{Prop}

\begin{Th}\label{T.1}
Let $(s,p) \in \R\times [1,+\infty)$.\\
For all $\varphi \in {\cal S}_e(\R^n)$ and $T\in E_{\alpha}^{s,p}$
we have  $\varphi .T \in E_{\alpha}^{s,p}$ and
\begin{eqnarray*}
 \Psi:&{\cal S}_e(\R^n)\times E_{\alpha}^{s,p}&\longrightarrow E_{\alpha}^{s,p}\\
 &(\varphi,T)&\longmapsto \varphi. T
 \end{eqnarray*}
 is bilinear continuous mapping.
\end{Th}
{\bf{Proof.}} Let $\varphi\in {\cal S}_e(\R^n)$ and $T\in
 E_{\alpha}^{s,p}$. Using relation (\ref{e1}), we have

$${\cal F}_{\alpha}(\varphi. T)(\xi)=c_{\alpha}^2\dint_{\R^n_+}
{\cal F}_{\alpha}(T)(x){\cal T}_{\xi}(\varphi)(x)d\mu_{\alpha}(x).$$
Using the following inequality (\cite{zuily}, p.186):
$$\Big(1+\|\xi\|^2\Big)^s \leq 2^{|s|}\Big(1+\|x\|^2\Big)^s\Big(1+\|x-\xi\|^2\Big)^{|s|},$$
we obtain that for all $\xi\in\R^n_+,$
$$(1+\|\xi\|^2)^{s}|{\cal F}_{\alpha}(\varphi. T)(\xi)|
\leq 2^{|s|}c_{\alpha}^2\dint_{\R^n_+} [(1+\|x\|^2)^{s}{\cal
F}_{\alpha}(T)(x)]\times [(1+(\|x-\xi\|)^2)^{|s|}|({\cal
T}_{\xi}{\cal F}_{\alpha}(\varphi))(x)|]d\mu_{\alpha}(x).$$ Since,
for all $x,\xi\in\R^n_+,$

$$(1+(\|x-\xi\|)^2)^{|s|}|({\cal T}_{\xi}{\cal F}_{\alpha}(\varphi))(x)|\leq
{\cal T}_{\xi}((1+\|t\|^2)^{|s|}|{\cal F}_{\alpha}(\varphi)|)(x).$$
 Then, we have
 \begin{eqnarray*}
(1+\|\xi\|^2)^{s}|{\cal F}_{\alpha}(\varphi. T)(\xi)| &\leq&
2^{|s|}c_{\alpha}^2\dint_{\R^n_+} [(1+\|x\|^2)^{s}{\cal
F}_{\alpha}(T)(x)]\times [{\cal T}_{\xi}((1+\|t\|^2)^{|s|}|{\cal
F}_{\alpha}(\varphi)|)(x)]d\mu_{\alpha}(x)\\
&=& 2^{|s|}c_{\alpha}^2\times [(1+\|x\|^2)^{s}{\cal
F}_{\alpha}(T)(x)]\ast_{\alpha} [(1+\|x\|^2)^{|s|}|{\cal
F}_{\alpha}(\varphi)|(\xi)].
\end{eqnarray*}

Therefore, again by (\ref{e1}) we obtain
$$\|(1+\|\xi\|^2)^{s}{\cal F}_{\alpha}(\varphi.T)\|_{L_{\alpha}^p}\leq
 2^{|s|}c_{\alpha}^2\times \|(1+\|x\|^2)^{|s|}{\cal F}_{\alpha}(T)
\|_{L_{\alpha}^p}.\|(1+\|x\|^2)^{s}{\cal
F}_{\alpha}(\varphi)\|_{L_{\alpha}^{1}}$$ i.e
\begin{equation}\label{e:I.3}
\|\varphi .T\|_{ E_{\alpha}^{s,p}}\leq 2^{|s|}c_{\alpha}\|T\|_{
E_{\alpha}^{s,p}} \|(1+\|x\|^2)^{s}{\cal
F}_{\alpha}(\varphi)\|_{L_{\alpha}^{1}}.
\end{equation}

On the other hand, we have:
\begin{eqnarray*}
\|(1+\|x\|^2)^{s}{\cal
F}_{\alpha}(\varphi)\|_{L_{\alpha}^{1}}&=&\dint_{\R^n_+}
(1+\|x\|^2)^{s}|{\cal
F}_{\alpha}(\varphi)(x)|d\mu_{\alpha}(x)\\
&\leq&
\Big(\dint_{\R^n_+}\frac{\prod_{i=1}^{n}x_i^{2\alpha_i+1}}{(1+\|x\|^2)^{s+|\alpha|
+\frac{3n}{2}}}dx_1...dx_n\Big).\sup_{x\in\R^n_+}(1+\|x\|^2)^{\nu}|{\cal
F}_{\alpha}(\varphi)(x)|,
\end{eqnarray*}
where $\nu\geq s+|\alpha|+\frac{3n}{2}$. Since the Fourier-Bessel
transform is continuous from ${\cal S}_e(\R^n)$ into itself, there
exist a positive integer $\nu'$ and a constant $c'\geq 0$, such that
 $$\sup_{x\in \R^n_+}(1+\|x\|^2)^{\nu}|{\cal
F}_{\alpha}(\varphi)(x)|\leq\sup_{x\in
\R^n_+;p,|\sigma|\leq\nu'}(1+\|x\|^2)^{p}|\partial^{\sigma}\varphi(x)|.$$

Then, $$\|\varphi. T\|_{ E_{\alpha}^{s,p}}\leq
c_{\alpha}c'\Big(\dint_{\R^n_+}
\frac{d\mu_{\alpha}(x)}{(1+\|x\|^2)^{s+|\alpha|+\frac{3n}{2}}}
\Big)\|T\|_{ E_{\alpha}^{s,p}}\rho_{\nu'}(\varphi).$$

\begin{Prop}\label{p:I.3}
\begin{enumerate}
\item Let $p\in [1,+\infty)\;and\; s,t \in \R$ such that $s\leq
t$, then $$E_{\alpha}^{t,p}\subset E_{\alpha}^{s,p}.$$ Moreover,
$\|T\|_{E_{\alpha}^{s,p}}\leq \|T\|_{E_{\alpha}^{t,p}}$ for all
$T\in E_{\alpha}^{t,p}$,  so that the embedding
$E_{\alpha}^{t,p}\hookrightarrow E_{\alpha}^{s,p}$ is continuous.
\item Let $(s,p) \in \R\times [1,+\infty)$, we have
$$(-\Delta_{\alpha})^k(E_{\alpha}^{s,p})\subset E_{\alpha}^{s-k,p},\;\;for \;all \;k\in \N.$$
Moreover, the operator
$(-\Delta_{\alpha})^k:E_{\alpha}^{s,p}\longrightarrow
E_{\alpha}^{s-k,p}$ is continuous and we have
$$\|(-\Delta_{\alpha})^k(T)\|_{E_{\alpha}^{s-k,p}}\leq \|T\|_{E_{\alpha}^{s,p}},\;\; T\in E_{\alpha}^{s,p}.$$
\end{enumerate}
\end{Prop}

{\bf{Application:}} Regularity of the solution of the equation
$P(-\Delta_{\alpha})T=u,$ where $u\in E_{\alpha}^{s,2}$ and $P$ an
even polynomial.

\begin{Prop}\label{p:I.4}
Let $s \in \R$, $u\in E_{\alpha}^{s,2}$ and $g \in
L_{\alpha}^2(\R^n_+)$. Assume that $P(-\Delta_{\alpha})g=u$, where P
is an even polynomial of degree 2m, then $g\in E_{\alpha}^{s+m,2}$.
\end{Prop}

{\bf{Proof.}} Using the relation (\ref{e:e}), we obtain:
$${\cal F}_{\alpha}(u)(\xi)={\cal F}_{\alpha}(P(-\Delta_{\alpha})g)(\xi)=P(\xi){\cal F}_{\alpha}(g)(\xi).$$
Since $u\in E_{\alpha}^{s,2}$, we deduce that
\begin{equation}\label{e:I.4}
\dint_{\R^n_+} (1+\|\xi\|^2)^{2s}|P(\xi)|^2 |{\cal
F}_{\alpha}(g)(\xi)|^2d\mu_{\alpha}(\xi)<\infty.
\end{equation}
On the other hand, we have
\begin{equation}\label{e:I.5}
(1+\|\xi\|^2)^{2(s+m)}|{\cal F}_{\alpha}(g)(\xi)|^2\sim |{\cal
F}_{\alpha}(g)(\xi)|^2 \quad (\|\xi\|\longrightarrow 0),
\end{equation}
and
\begin{equation}\label{e:I.6}
(1+\|\xi\|^2)^{2(s+m)}|{\cal F}_{\alpha}(g)(\xi)|^2\sim
(1+\|\xi\|^2)^{2s}|P(\xi)|^2|{\cal F}_{\alpha}(g)(\xi)|^2 \quad
(\|\xi\|\longrightarrow \infty).
\end{equation}
Since $g\in L_{\alpha}^2(\R^n_+)$, using the relation (\ref{e:I.5})
we obtain
$$\dint_{0<\|\xi\|\leq 1}(1+\|\xi\|^2)^{2(s+m)}|{\cal F}_{\alpha}(g)(\xi)|^2d\mu_{\alpha}(\xi)<\infty.$$
From relations (\ref{e:I.4}) and (\ref{e:I.6}) we deduce that
$$\dint_{\|\xi\|>1}(1+\|\xi\|^2)^{2(s+m)}|{\cal F}_{\alpha}(g)(\xi)|^2d\mu_{\alpha}(\xi)<\infty.$$
This proves that $g\in E_{\alpha}^{s+m,2}$.

%%%%%%%%%%%%%%%%%%%%%%%%%%%%%%%%%%%%%%%%%%%%%%%%%%%%%%%%%%%%%%%%%%%%%%%%%%%%%%%%
\section{   Special case of the Sobolev type spaces $E_{\alpha}^{s,2}$}
%%%%%%%%%%%%%%%%%%%%%%%%%%%%%%%%%%%%%%%%%%%%%%%%%%%%%%%%%%%%%%%%%%%%%%%%%%%%%%%%

This section is devoted to the study of the Hilbert spaces
$E_{\alpha}^{s,2}$ which will be denoted in the sequel by
$H_{\alpha}^s$.

\begin{Prop}\label{p:II.1}
Let $m\in \N$. Then for all $s>\frac{1}{2}(|\alpha|+n)+m,$
$$H_{\alpha}^s \subset {\cal C}_e^m(\R^n).$$
\end{Prop}
To prove this proposition we need the following lemma.
\begin{Lem}\label{l1}
Let $f:\R^n_+\longrightarrow \C$ such that for all $k\in
\left\{0,...,m\right\}$, $\|x\|^{2k}f \in L_{\alpha}^1(\R^n_+)$ then
${\cal F}_{\alpha}(f) \in {\cal C}_e^m(\R^n).$
\end{Lem}

{\bf{Proof of Proposition \ref{p:II.1}.}} Let $f\in H_{\alpha}^s$,
by relation (\ref{e:e}) and by lemma \ref{l1} it suffices to show
that
 $$\|\xi\|^{2k}{\cal F}_{\alpha}(f)\in L_{\alpha}^1(\R^n_+),\;\;for\;all\;k\in\left\{0,...,m\right\}.$$
Let $k\in\left\{0,...,m\right\}$ we have:  $\|\xi\|^{2k}{\cal
F}_{\alpha}(f)(\xi)={\cal
F}_{\alpha}((-\Delta_{\alpha})^kf)(\xi)$.\\
Then for $f\in H_{\alpha}^s$ and $k\in\left\{0,...,m\right\}$
$$\dint_{\R^n_+} (1+\|\xi\|^2)^{2(s-k)}|{\cal F}_{\alpha}((-\Delta_{\alpha})^kf)(\xi)|^2d\mu_{\alpha}(\xi)<\infty.$$
On the other hand, by applying H\"{o}lder inequality, for all
$k\in\left\{0,...,m\right\}$

$$\dint_{\R^n_+}|{\cal
F}_{\alpha}((-\Delta_{\alpha})^kf)(\xi)|^2d\mu_{\alpha}(\xi)
 \leq
 \Big( \dint_{\R^n_+} (1+\|\xi\|^2)^{-2(s-k)}d\mu_{\alpha}(\xi) \Big)^{1/2}$$
 $$\times  \Big(\dint_{\R^n_+} (1+\|\xi\|^2)^{2(s-k)} |{\cal
F}_{\alpha}((-\Delta_{\alpha})^kf)(\xi)|^2d\mu_{\alpha}(\xi)
\Big)^{1/2}.$$ Note that the latter integral converges since
$s>\frac{1}{2}(|\alpha|+n)+m $,
$$\dint_{\R^n_+}(1+\|\xi\|^2)^{-2(s-k)}d\mu_{\alpha}(\xi)<\infty.$$
Therefore, for all $k\in\left\{0,...,m\right\}$,
$\dint_{\R^n_+}|{\cal
F}_{\alpha}((-\Delta_{\alpha})^kf)(\xi)|^2d\mu_{\alpha}(\xi)<\infty$.\\
This completes the proof of the proposition.\\
\\
{\bf{Notation.}} For all $s \in \R$, we denoted by
$(H_{\alpha}^s)^*$ the topological dual of
$(H_{\alpha}^s;\|.\|_{H_{\alpha}^s})$ and by
$\|.\|_{(H_{\alpha}^s)^*}$ the usual uniform norm on
$(H_{\alpha}^s)^*$ given by:
$$\|L\|_{(H_{\alpha}^s)^*}=\dsup_{T\in H_{\alpha}^s, \|T\|_{H_{\alpha}^s}\leq 1}\|L(T)\|_{H_{\alpha}^s}.$$

\begin{Th}\label{T:I.1}
Let $s \in \R$, then the following holds:
\begin{enumerate}
\item Every tempered distribution $T\in H_{\alpha}^{-s}$
extends uniquely in a continuous linear form $L_T$ on
$(H_{\alpha}^s,\|.\|_{H_{\alpha}^s})$.
\item The map
$$\chi:\Big(H_{\alpha}^{-s},\|.\|_{H_{\alpha}^{-s}}\Big)\longrightarrow
\Big((H_{\alpha}^s)^*,\|.\|_{(H_{\alpha}^s)^*}\Big)$$
$$T\longmapsto L_T.$$ is an isometric isomorphism.
\end{enumerate}
\end{Th}

{\bf{Proof.}} 1. By inversion formula, theorem \ref{IF}, for the
Fourier-Bessel transform ${\cal F}_{\alpha}$ we have
$$\forall \varphi \in {\cal S}_e(\R^n);\;\;\varphi=c_{\alpha}^2
{\cal F}_{\alpha}({\cal F}_{\alpha}(\varphi)).$$ Since, for all
$T\in H_{\alpha}^{-s}$ and $\varphi \in {\cal S}_e(\R^n)$,
\begin{eqnarray*}
<T,\varphi>&=&c_{\alpha}^2<{\cal F}_{\alpha}(T),{\cal F}_{\alpha}(\varphi)>\\
&=&c_{\alpha}^2\dint_{\R^n_+}
{\cal F}_{\alpha}(T)(\xi){\cal F}_{\alpha}(\varphi)(\xi)d\mu_{\alpha}(\xi)\\
&=&c_{\alpha}^2\dint_{\R^n_+} (1+\|\xi\|^2)^{s}{\cal
F}_{\alpha}(\varphi)(\xi) \Big[(1+\|\xi\|^2)^{-s} {\cal
F}_{\alpha}(T)(\xi) \Big]d\mu_{\alpha}(\xi).
\end{eqnarray*}
Using H\"{o}lder inequality, we obtain
$$\forall \varphi \in {\cal S}_e(\R^n);\;\;|<T,\varphi>|\leq
\|\varphi\|_{H_{\alpha}^{s}}.\|T\|_{H_{\alpha}^{-s}}.$$ By the
density of ${\cal S}_e(\R^n)$ in $H_{\alpha}^{s}$ and the Hahn
Banach theorem \cite{hb} we deduce that T extends uniquely to a
continuous linear form
 $L_T:(H_{\alpha}^{s},\|.\|_{H_{\alpha}^{s}})\longrightarrow
\C$ such that
\begin{equation}
\|L_T\|_{(H_{\alpha}^{s})*}\leq \|T\|_{H_{\alpha}^{-s}}.
\end{equation}
2. The linearity and injectivity of $\chi$ can be deduced directly
from the uniqueness of the extension of each $T\in H_{\alpha}^{-s}$
in a continuous linear
form $L_T=\chi(T)\in (H_{\alpha}^{-s})^*$.\\
It suffices to prove that $\chi$ is surjective. Let
$L:H_{\alpha}^{s}\longrightarrow \C $ be a continuous linear form.\\
By proposition \ref{p:I.2} and the Riesz theorem \cite{hb}, there
exists $S\in H_{\alpha}^{s}$ such that
$$\|S\|_{H_{\alpha}^{s}}=\|L\|_{(H_{\alpha}^{s})*}\quad\;and\quad\;
\forall u\in H_{\alpha}^{s},\;L(u)=(u,S)_{H_{\alpha}^{s}}$$

$\forall \varphi\in H_{\alpha}^{s},\;L(\varphi)=<T,\varphi>$
 with $T$ is the tempered distribution defined by:
$$T={\cal F}^{-1}_{\alpha}((1+\|\xi\|^2)^{2s}\overline{{\cal F}_{\alpha}(S)}).$$
This equality is equivalent to:
\begin{equation}\label{e.a}
(1+\|\xi\|^2)^{-s}{\cal
F}_{\alpha}(T)=(1+\|\xi\|^2)^{s}\overline{{\cal F}_{\alpha}(S)},
\end{equation}
since $T \in H_{\alpha}^{-s}$, then $L=L_T$.\\ It remains to prove
that $\chi$ is an isometry, from (\ref{e.a}) we have
$$\|S\|_{H_{\alpha}^{s}}=\|T\|_{H_{\alpha}^{-s}}$$
since
$\|S\|_{H_{\alpha}^{s}}=\|L\|_{(H_{\alpha}^{s})*}\;\;and\;\;T=\chi^{-1}(L)$,
therefore
$$\|L\|_{(H_{\alpha}^{s})*}=\|\chi^{-1}(L)\|_{H_{\alpha}^{s}}.$$
This leads to the result.\\
In what follows we give a new characterization of $H_{\alpha}^{-s}$
for $s\in N$.

\begin{Th}\label{T:I.2}
For all $m\in \N$, the Sobolev type space $H_{\alpha}^{-m}$ is
generated by the set
$$\left\{(-\Delta_{\alpha})^kg,\;for \;all \;g\in
L_{\alpha}^2(\R^n_+)\;\;and\;k\in \left\{0,...,m\right\}\right\}.$$
 Moreover, for all $T\in
H_{\alpha}^{-m}$ there exists $g\in L_{\alpha}^2(\R^n_+)$ such that
$$T=(1-\Delta_{\alpha})^mg=\sum_{k=0}^m C_m^k(-\Delta_{\alpha})^kg. $$
\end{Th}

{\bf{Proof.}} Since
 $L_{\alpha}^2(\R^n_+)=H_{\alpha}^0(\R^n_+)=E_{\alpha}^{0,2}(\R^n_+)$
 then $(-\Delta_{\alpha})^kg\in
 H_{\alpha}^{-k}=E_{\alpha}^{-k,2}$.\\
For all $g\in L_{\alpha}^2(\R^n_+)$ and $k\in
\left\{0,...,m\right\}$, by proposition \ref{p:I.3} we have
$$\left\{(-\Delta_{\alpha})^kg,\;for\;all\;\;g\in
L_{\alpha}^2(\R^n_+)\;\;and\;k\in
\left\{0,...,m\right\}\right\}\subset H_{\alpha}^{-m}.$$  Let now $T
\in H_{\alpha}^{-m}$ then by Plancherel theorem, there exists $g\in
L_{\alpha}^2(\R^n_+)$ such that
$$(1+\|\xi\|^2)^{-m}{\cal F}_{\alpha}(T)={\cal F}_{\alpha}(g).$$
This gives the result since the Fourier-Bessel transform is
injective.

\begin{Prop}\label{p:II.2}
The following holds for all $\varphi \in {\cal S}_e(\R^n)$ and $s,t
\in \R$ such that $t<s$, then
$$H_{\alpha}^s\longrightarrow H_{\alpha}^t$$
$$T\longmapsto \varphi\;T$$
is a compact operator.
\end{Prop}

{\bf{Proof.}} Let $(T_n)_{n\in \N}$ be a sequence of $H_{\alpha}^s$
such that $\|T_n\|_{H_{\alpha}^s}\leq 1$. By the Alaoglu theorem
(\cite{hb},p.42) we can extract a subsequence $(T_{n_{k}})_{k\in
\N}$ weakly convergent to $T$ in $H_{\alpha}^s$. Then
\begin{equation}\label{e:II.3}
\forall S\in H_{\alpha}^s,\; \dlim_{k\longrightarrow
\infty}(T_{n_{k}},S)_{H_{\alpha}^s} =(T,S)_{H_{\alpha}^s}.
\end{equation}
Prove that the sequence $(\varphi\;T_{n_{k}})_{k\in \N}$ converges
to $(\varphi\;T)$ in $H_{\alpha}^t$ is equivalent to show that:
\begin{equation}
\dlim_{k\longrightarrow\infty}\|\varphi\;v_k\|^2_{H_{\alpha}^t}=0,
\end{equation}
where $v_k=T_{n_{k}}-T$ for all $k\in \N$, then for all $R>0$ we
have
\begin{equation}\label{e.1}
\|\varphi\;v_k\|^2_{H_{\alpha}^t}\leq c_{\alpha}^2
 \Big(\dint_{0<\|\xi\|\leq R} (1+\|\xi\|^2)^{2t}|{\cal
F}_{\alpha}(\varphi\;v_k)|^2d\mu_{\alpha}(\xi) +\dint_{\|\xi\|>
R}(1+\|\xi\|^2)^{2t}|{\cal
F}_{\alpha}(\varphi\;v_k)|^2d\mu_{\alpha}(\xi)\Big).
\end{equation}
Since
\begin{eqnarray*}
\dint_{\|\xi\|> R}(1+\|\xi\|^2)^{2t}|{\cal
F}_{\alpha}(\varphi\;v_k)|^2d\mu_{\alpha}(\xi)
&\leq&\frac{\|\varphi\;v_k\|^2_{H_{\alpha}^s}}{(1+R^2)^{2(s-t)}}.
\end{eqnarray*}
On the other hand, by theorem \ref{T.1} there exist a constant $C>0$
and a positive integer $\nu$ such that
$$\forall k\in \N,\;\|\varphi\; v_k\|^2_{H_{\alpha}^s}\leq
C(\rho_\nu(\varphi))^2\|v_k\|^2_{H_{\alpha}^s}.$$ So that
$$\forall k\in\N,\;\|\varphi\; v_k\|^2_{H_{\alpha}^s}\leq
C(\rho_\nu(\varphi))^2(1+\|T\|_{H_{\alpha}^s})^2.$$ Thus

$$\;\;\forall k\in \N,\;\; c_{\alpha}^2\dint_{\|\xi\|>
R}(1+\|\xi\|^2)^{2t} |{\cal F}_{\alpha}(\varphi\;
v_k)|^2d\mu_{\alpha}(\xi)\leq
\frac{c_{\alpha}^2C(\rho_\nu(\varphi))^2(1+\|T\|_{H_{\alpha}^s})^2}{(1+R^2)^{2(s-t)}}.$$
Hence for all $\varepsilon >0$ and taking $R$ large we have
$$\frac{c_{\alpha}^2C(\rho_\nu(\varphi))^2(1+\|T\|_{H_{\alpha}^s})^2}{(1+R^2)^{2(s-t)}}
<\frac{\varepsilon}{2}.$$ And relation (\ref{e.1}) gives
\begin{equation}\label{e:II.4}
\forall k\in \N,\;\;\|\varphi\; v_k\|^2_{H_{\alpha}^t}\leq
 \Big(c_{\alpha}^2\dint_{0<\|\xi\|\leq R}(1+\|\xi\|^2)^{2t}|{\cal
F}_{\alpha}(\varphi\; v_k)|^2d\mu_{\alpha}(\xi)
\Big)+\frac{\varepsilon}{2}.
\end{equation}
We show that
\begin{equation}
\dlim_{k\longrightarrow \infty}
\Big(c_{\alpha}^2\dint_{0<\|\xi\|\leq R}(1+\|\xi\|^2)^{2t}|{\cal
F}_{\alpha}(\varphi\;v_k)|^2d\mu_{\alpha}(\xi) \Big)=0.
\end{equation}
From relation (\ref{e1}) we have
\begin{equation}\label{e.2}
\forall k\in \N,\quad{\cal F}_{\alpha}(\varphi\;v_k)=c_{\alpha}^2
 \Big({\cal F}_{\alpha}(v_k)\ast_{\alpha} {\cal F}_{\alpha}(\varphi) \Big).
\end{equation}
In the sense of distribution
\begin{eqnarray*}
{\cal F}_{\alpha}(\varphi\;v_k)(\xi)&=&c_{\alpha}^2 \langle{\cal
F}_{\alpha}(v_k),({T}_{\xi}{\cal
F}_{\alpha}(\varphi))\rangle\\
&=&c_{\alpha}^2 \langle (1+\|x\|^2)^{s}{\cal
F}_{\alpha}(v_k),(1+\|x\|^2)^{-s}({T}_{\xi}{\cal
F}_{\alpha}(\varphi))\rangle.
\end{eqnarray*}

Using the fact that
\begin{equation}\label{e:II.5}
{\cal F}_{\alpha}(\varphi\;v_k)(\xi)= \Big(v_k,{\cal
F}_{\alpha}^{-1}\Big[(1+\|x\|^2)^{-2s}({T}_{\xi}\overline{{\cal
F}_{\alpha}(\varphi)})\Big] \Big)_{H_{\alpha}^s}.
\end{equation}
Since the sequence $(v_k)_k$ is weakly convergent to $0$ in
$H_{\alpha}^s$, we deduce that
\begin{equation}\label{e:II.6}
\forall \xi\in\R^n_+,\;\dlim_{k\longrightarrow \infty}{\cal
F}_{\alpha}(\varphi\;v_k)(\xi)=0.
\end{equation}
From the relation (\ref{e:II.5}) and by Cauchy-Schwartz inequality
we deduce that
\begin{equation}\label{e:II.7}
\forall k\in \N ,\;\;\xi\in \R^n_+,\;\;\; |{\cal
F}_{\alpha}(\varphi\;v_k)(\xi)| \leq C_T
\|(1+\|x\|^2)^{-2s}({T}_{\xi} \overline{{\cal
F}_{\alpha}(\varphi)})\|_{H_{\alpha}^s}
\end{equation}
where $C_T=(1+\|T\|_{H_{\alpha}^s})$. Since
$$\|(1+\|x\|^2)^{-2s}({T}_{\xi} \overline{{\cal
F}_{\alpha}(\varphi)})\|^2_{H_{\alpha}^s}\leq\|\varphi\|^2_{L^2_{\alpha}}.$$
By the Plancherel theorem we have
\begin{equation}\label{e:II.8}
\forall k\in \N ,\;\xi\in \R^n_+, |{\cal
F}_{\alpha}(\varphi\;v_k)(\xi)|^2\leq C_T
\|\varphi\|^2_{L^2_{\alpha}}.
\end{equation}
By the dominated convergence theorem and relations (\ref{e:II.6})
and (\ref{e:II.8}) we deduce
$$\dlim_{k\longrightarrow \infty} \Big(c_{\alpha}^2\dint_{0<\|\xi\|\leq R}(1+\|\xi\|^2)^{2t}|{\cal
F}_{\alpha}(\varphi\;v_k)|^2d\mu_{\alpha}(\xi) \Big)=0.$$ Thanks to
relation (\ref{e:II.4}) we have: for all $\varepsilon
>0$
$$\|\varphi\;v_k\|^2_{H_{\alpha}^t}\leq \varepsilon.$$
Then the result is proved.\\

{\bf{Notation.}} Let $K$ be a compact set of $\R^n$, we denoted by
$H_{\alpha,K}^s,\;s\in \R$, the subspace of $H_{\alpha}^s$ defined
by:
$$H_{\alpha,K}^s=\left\{T\in H_{\alpha}^s/ supp(T) \subset K\right\}$$
with $supp(T)$ is the support of $T$.\\

{\bf{Remark.}} For all $s\geq 0$ and all compact $K\subset \R^n_+$,
we have: $H_{\alpha,K}^s\subset L_{\alpha,K}^2(\R^n_+).$\\
with $L_{\alpha,K}^2(\R^n_+)=\left\{f\in
L_{\alpha}^2(\R^n_+)/\;\;f(x)=0,\;\;x\in\R^n_+\setminus K\right\}$.

\begin{Th}[ Reillich theorem].
Let $s,t\in \R$ such that $\;t<s$, for all compact $K\subset \R^n$
the embedding $H_{\alpha,K}^s\hookrightarrow H _{\alpha,K}^t$ is a
compact operator.
\end{Th}

{\bf{Proof.}} Let $\widetilde{K}=K\cup(-K)$ be a compact of $\R^n$,
$V$ a relatively compact neighborhood of $K$, $\widetilde{V}$ a
relatively compact neighborhood of $\widetilde{K}$. By the Urysohn
theorem (\cite{gf},p.237) there exists a $\varphi \in {\cal
D}_e(\R^n)$ such
that $\varphi\equiv 1$ in $\widetilde{V}$.\\
Put: $$\quad\quad\Psi(x)=\frac{\varphi(x)+\varphi(-x)}{2},\;\; for\;
all\; x\in \R^n_+.$$ Therefore $\quad\quad\Psi .T=T\quad for\; all\;
T\in
 H_{\alpha,K}^s\;(resp. \;\;H_{\alpha,K}^t)$.\\
By proposition \ref{p:II.2}, we deduce the result.

\begin{Cor}\label{c:II.1}
Let $K\subset \R^n$ be a compact set, for all $s\geq 0$ there exists
a constant $C>0$ such that for all $T\in H_{\alpha,K}^s$, we obtain
$$\frac{1}{C}\|T\|_{H_{\alpha}^s}\leq
c_{\alpha} \Big(\dint_{\R^n_+}\|\xi\|^{4s} |{\cal
F}_{\alpha}(T)(\xi)|^2d\mu_{\alpha}(\xi) \Big)^{1/2}\leq
C.\|T\|_{H_{\alpha}^s}.$$
\end{Cor}

{\bf{Proof.}} Suppose that a such constant does not exists. Then,
for all $k\in \{1,2,3...\}$ there exists $T_k\in H_{\alpha,K}^s$
such that
$$\frac{1}{k}\|T_k\|_{H_{\alpha}^s}>c_{\alpha} \Big(\dint_{\R^n_+}\|\xi\|^{4s}
|{\cal F}_{\alpha}(T)(\xi)|^2d\mu_{\alpha}(\xi) \Big)^{1/2}.$$ Put
$\quad\widetilde{T_k}=\frac{1}{\|T_k\|_{H_{\alpha}^s}}T_k$, for all
$k\in\{1,2,...\}$, therefore
\begin{equation}\label{e:II.10}
 \|\widetilde{T_k}\|_{H_{\alpha}^s}=1,\quad for\;all\;k\in\{1,2,3...\}
\end{equation}
and we have:
$$ \Big(\dint_{\R^n_+}\|\xi\|^{4s}
|{\cal F}_{\alpha}(\widetilde{T_k})(\xi)|^2d\mu_{\alpha}(\xi)
\Big)^{1/2}\leq \frac{c_{\alpha}^{-1}}{k},\quad
for\;all\;k\in\{1,2,3...\}.$$ So it holds
\begin{equation}\label{e:II.11}
\dlim_{k\longrightarrow \infty} \Big(\dint_{\R^n_+}\|\xi\|^{4s}
|{\cal F}_{\alpha}(\widetilde{T_k})(\xi)|^2d\mu_{\alpha}(\xi)
\Big)=0.
\end{equation}
From Reillich theorem we deduce that the sequence
$(\widetilde{T_k})_{k\geq 1}$ admits a subsequence
$(\widetilde{{T_k}_p})_{p\in \N}$ convergent to $\widetilde{T}$ in $
L_{\alpha,K}^2(\R^n_+)$. By H\"{o}lder inequality we have

$$\|\widetilde{{T_k}_p}-\widetilde{T}\|_{L_{\alpha,}^1}\leq
C_{K}\|\widetilde{{T_k}_p}-\widetilde{T}\|_{L_{\alpha}^2}.$$ Since
$$\| {\cal F}_{\alpha}(\widetilde{{T_k}_p})-{\cal
F}_{\alpha}(\widetilde{T})\|_{L_{\alpha}^\infty}\leq
\|\widetilde{{T_k}_p}-\widetilde{T}\|_{L_{\alpha}^1}$$

which implies that $\dlim_{k\rightarrow\infty}\| {\cal
F}_{\alpha}(\widetilde{{T_k}_p})-{\cal
F}_{\alpha}(\widetilde{T})\|_{L_{\alpha}^\infty}=0$. So that
$$\forall\;R>0,\;\;\dlim_{p\rightarrow \infty}\dint_{0<\|\xi\|\leq R}\|\xi\|^{4s}|{\cal F}_{\alpha}(\widetilde{{T_k}_p})(\xi)|^2d\mu_{\alpha}
(\xi)=\dint_{0<\|\xi\|\leq R}\|\xi\|^{4s}|{\cal
F}_{\alpha}(\widetilde{T})(\xi)|^2d\mu_{\alpha} (\xi).$$ Then we
tend R to $+\infty$, we obtain
$$\dint_{\R^n_+}\|\xi\|^{4s}|{\cal F}_{\alpha}(\widetilde{T})(\xi)|^2d\mu_{\alpha}
(\xi)=0.$$ and hence $\widetilde{T}=0$, % (because the Fourier-Bessel transform is injective ).
then the sequence $(\widetilde{{T_k}_p})_{p\in \N}$ converges to $0$
in $L_{\alpha}^2(\R^n_+)$.\\
On the other hand, the function $\quad\xi\longmapsto
\frac{(1+\|\xi\|^2)^{2s}}{1+\|\xi\|^{4s}}$ is continuous and bounded
on $\R^n_+$.\\
Thus, there exists $\lambda>0$ such that:
$$\forall \xi\in \R^n_+,\;(1+\|\xi\|^2)^{2s}\leq \lambda (1+\|\xi\|^{4s}).$$
Hence, for all $p\in \N$, we have:
$$\|\widetilde{{T_k}_p}\|^2_{H_{\alpha}^s}\leq c_{\alpha}^2\lambda
\dint_{\R^n_+}|{\cal
F}_{\alpha}(\widetilde{{T_k}_p})(\xi)|^2d\mu_{\alpha}
(\xi)+c_{\alpha}^2\lambda \dint_{\R^n_+}\|\xi\|^{4s}|{\cal
F}_{\alpha}(\widetilde{{T_k}_p})(\xi)|^2d\mu_{\alpha} (\xi).$$

From the relation (\ref{e:II.10}) and by Plancherel theorem we
obtain
\begin{equation}\label{e:II.13}
 1\leq \lambda
\|\widetilde{{T_k}_p}\|_{L_{\alpha}^2}^2+c_{\alpha}^2\lambda
\dint_{\R^n_+}\|\xi\|^{4s}|{\cal
F}_{\alpha}(\widetilde{{T_k}_p})(\xi)|^2d\mu_{\alpha}
(\xi),\;\;\forall p\in \N.
\end{equation}
which we obtain the absurdity by tending $p\rightarrow\infty$. Then,
we have
\begin{equation}\label{e:II.15}
\forall T \in H_{\alpha,K}^s,\;c_{\alpha}
\Big(\dint_{\R^n_+}\|\xi\|^{4s} |{\cal
F}_{\alpha}(T)(\xi)|^2d\mu_{\alpha} (\xi) \Big)^{1/2}\leq
\|T\|_{H_{\alpha}^s}.
\end{equation}
So, we deduce that
\begin{eqnarray*}
\forall T \in
H_{\alpha,K}^s,\quad\;\frac{1}{C+1}\|T\|_{H_{\alpha}^s}&\leq&c_{\alpha}
\Big(\dint_{\R^n_+}\|\xi\|^{4s} |{\cal
F}_{\alpha}(T)(\xi)|^2d\mu_{\alpha}
(\xi) \Big)^{1/2}\\
&\leq& (1+C) \|T\|_{H_{\alpha}^s}.
\end{eqnarray*}

\begin{Th}[Poincar\'e inequality].
Let $0\leq t\leq s$, then there exists a constant $C>0$ such that
$\forall \varepsilon \in\R_+;$
$$\|T\|_{H_{\alpha}^t}
\leq C.\varepsilon^{2(s-t)} \|T\|_{H_{\alpha}^s},\quad for\; all
\;\;T\in H_{\alpha,\varepsilon}^s.$$

\end{Th}

{\bf{Proof.}}  From the corollary \ref{c:II.1}, there exist two
positive constants $C_1$ and $C_2$ such that:
\begin{equation}\label{e:II.1.1}
\frac{1}{C_1}\|T\|_{H_{\alpha}^s}\leq c_{\alpha}
\Big(\dint_{\R^n_+}\|\xi\|^{4s} |{\cal
F}_{\alpha}(T)(\xi)|^2d\mu_{\alpha}(\xi) \Big)^{1/2}\leq
C_1.\|T\|_{H_{\alpha}^s},
\end{equation}
and
\begin{equation}\label{e:II.1.2}
\frac{1}{C_2}\|T\|_{H_{\alpha}^t}\leq c_{\alpha}
\Big(\dint_{\R^n_+}\|\xi\|^{4t} |{\cal
F}_{\alpha}(T)(\xi)|^2d\mu_{\alpha}(\xi) \Big)^{1/2}\leq
C_2.\|T\|_{H_{\alpha}^t}.
\end{equation}
On the other hand, we have:
\begin{equation}\label{e:II.1.3}
T\in H_{\alpha}^s,\quad \|T\|_{H_{\alpha}^t}\leq \|T\|_{
H_{\alpha}^s}.
\end{equation}
 From (\ref{e:II.1.1}),(\ref{e:II.1.2}) and (\ref{e:II.1.3}), we deduce that: for all $T\in H_{\alpha,1}^s$
\begin{equation}\label{e:II.16}
 \Big(\dint_{\R^n_+}\|\xi\|^{4t} |{\cal F}_{\alpha}(T)(\xi)|^2d\mu_{\alpha}(\xi) \Big)^{1/2}\leq C_1C_2 \Big(\dint_{\R^n_+}\|\xi\|^{4s} |{\cal
F}_{\alpha}(T)(\xi)|^2d\mu_{\alpha}(\xi) \Big)^{1/2}.
\end{equation}
Let $\varepsilon \in ]0,1[$ and $T\in H_{\alpha,\varepsilon}^s$, we
have the function defined by:
 $T_\varepsilon(x)=T(\varepsilon x_1,...,\varepsilon x_n),\;\;x\in
\R^n_+$. Then
\begin{eqnarray*}
\forall \xi\in\R^n_+,\quad{\cal F}_{\alpha}(T_\varepsilon)(\xi)
&=&\frac{1}{\varepsilon^{2|\alpha|+n}}{\cal
F}_{\alpha}(T)(\frac{\xi_1}{\varepsilon},...,\frac{\xi_n}{\varepsilon}).
\end{eqnarray*}

By applying the relation (\ref{e:II.16}) we have
$$ \Big(\dint_{\R^n_+}\|\xi\|^{4t}
|{\cal
F}_{\alpha}(T)(\frac{\xi_1}{\varepsilon},...,\frac{\xi_n}{\varepsilon})|^2d\mu_{\alpha}(\xi)
\Big)^{1/2} \leq C_1C_2 \Big(\dint_{\R^n_+}\|\xi\|^{4s} |{\cal
F}_{\alpha}(T)(\frac{\xi_1}{\varepsilon},...,\frac{\xi_n}{\varepsilon})|^2d\mu_{\alpha}(\xi)
\Big)^{1/2}.$$
 By the following change
 $\eta_i=\frac{\xi_i}{\varepsilon}$ for $i=1,...,n$;  we deduce that
\begin{equation}\label{e:II.17}
\varepsilon^{2t} \Big(\dint_{\R^n_+}\|\eta\|^{4t} |{\cal
F}_{\alpha}(T)(\eta)|^2d\mu_{\alpha}(\eta) \Big)^{1/2} \leq
C_1C_2\varepsilon^{2s} \Big(\dint_{\R^n_+}\|\eta\|^{4s}|{\cal
F}_{\alpha}(T)(\eta)|^2d\mu_{\alpha}(\eta) \Big)^{1/2}.
\end{equation}
From corollary \ref{c:II.1} there exist two positive constants $C_3$
and $C_4$ such that for all $T \in H_{\alpha,\varepsilon}^s$, we
have
\begin{equation}\label{e:II.18}
\|T\|_{H_{\alpha}^t}\leq C_3  \Big(\dint_{\R^n_+}\|\eta\|^{4t}
|{\cal F}_{\alpha}(T)(\eta)|^2d\mu_{\alpha}(\eta) \Big)^{1/2}
\end{equation}
and
\begin{equation}\label{e:II.19}
 \Big(\dint_{\R^n_+}\|\eta\|^{4s} |{\cal F}_{\alpha}(T)(\eta)|^2d\mu_{\alpha}(\eta) \Big)^{1/2}\leq C_4
 \|T\|_{H_{\alpha}^s}.
\end{equation}
 From relations (\ref{e:II.17}),(\ref{e:II.18})
and (\ref{e:II.19}) we deduce that: $\|T\|_{H_{\alpha}^t}\leq
C_1C_2C_3C_4 \varepsilon^{2(s-t)}\|T\|_{H_{\alpha}^s}$.\\

{\bf{Application.}} Regularity of solutions of the equation (E):
$(k^2-\Delta_{\alpha})u=f$.\\
For $f\in {\cal S}'_e(\R^n)$ and $k$ a non-zero constant, the
equation (E) has a unique tempered solution $u$, such that:
$${\cal F}_{\alpha}((k^2-\Delta_{\alpha})u)(x)=(k^2+\|x\|^2){\cal F}_{\alpha}(u)(x)={\cal F}_{\alpha}(f)(x);\quad x\in \R^n_+.$$
This is equivalent to
$${\cal F}_{\alpha}(u)=(k^2+\|x\|^2)^{-1}{\cal F}_{\alpha}(f).$$
Now, we suppose that $f \in H_{\alpha}^s$. Let $u\in {\cal
S}'_e(\R^n)$ the unique solution of (E), then
\begin{eqnarray*}
(1+\|x\|^2)^{s+1}{\cal
F}_{\alpha}(u)(x)&=&(\frac{1+\|x\|^2}{k^2+\|x\|^2})(1+\|x\|^2)^s(k^2+\|x\|^2){\cal F}_{\alpha}(u)(x)\\
&=&(\frac{1+\|x\|^2}{k^2+\|x\|^2})(1+\|x\|^2)^s {\cal
F}_{\alpha}(f)(x).
\end{eqnarray*}
Since $x\longmapsto \frac{1+\|x\|^2}{k^2+\|x\|^2}$ is a bounded
function, then $(1+\|x\|^2)^{s+1}{\cal F}_{\alpha}(u) \in
L_{\alpha}^2(\R^n_+)$.\\
So, for all $f \in H_{\alpha}^s$, the unique solution of (E) belongs
to $H_{\alpha}^{s+1}$.

\end{document}